\documentclass[11pt,oneside]{amsart}

\usepackage[draft=false]{hyperref}
\usepackage{amsmath,amsthm,amssymb,amsfonts}
\usepackage{etoolbox}
\usepackage{dsfont}
\usepackage{pifont}
\usepackage{mathrsfs}
\usepackage{graphicx}
\usepackage[shortlabels]{enumitem}
\usepackage[all]{xy}
\usepackage{geometry}
\usepackage{amsfonts}
\usepackage{fancyhdr}
\usepackage{tikz}
\usepackage{tikz-cd}
\usetikzlibrary{arrows,shapes,automata,backgrounds,petri}
\usepackage[section]{algorithm}
\usepackage{algorithmic}
\usepackage{lscape}
\usepackage{quiver}

\theoremstyle{plain}
\newtheorem{prop}{Proposition}[section]
\newtheorem{prop-def}[prop]{Proposition-Definition}
\newtheorem{lem}[prop]{Lemma}
\newtheorem{theo}[prop]{Theorem}
\newtheorem{cor}[prop]{Corollary}

\newtheorem*{theo*}{Theorem}
\newtheorem*{cor*}{Corollary}
\newtheorem*{prop*}{Proposition}
\newtheorem*{fact*}{Fact}

\theoremstyle{definition}
\newtheorem{definition}[prop]{Definition}

\theoremstyle{remark}
\newtheorem{rem}[prop]{Remark}

 \normalbaselines
\numberwithin{equation}{section}

\makeatletter
\def\thmhead@plain#1#2#3{%
  \thmname{#1}\thmnumber{\@ifnotempty{#1}{ }\@upn{#2}}%
  \thmnote{ {\the\thm@notefont#3}}}
\let\thmhead\thmhead@plain
\makeatother

\usepackage[backend=biber,style=alphabetic,sorting=nty,maxnames=5]{biblatex}
\newcommand{\C}{\mathbb{C}}

\renewcommand{\P}{\mathbb{P}}

\DeclareMathOperator\GL{GL}
\DeclareMathOperator\PGL{PGL}

\DeclareMathOperator\Id{Id}

\DeclareMathOperator\Stab{Stab}

\newcommand{\xdownarrow}[1]{%
  {\left\downarrow\vbox to #1{}\right.\kern-\nulldelimiterspace}
}

\renewcommand{\matrix}[1]{\begin{pmatrix} #1\end{pmatrix}}

\newcommand{\muls}[1]{ \left\{\kern-0.6em \left\{ #1\right\}\kern-0.6em \right\} }%Multisets

\newcommand{\hhat}[1]{\widehat{#1}}
\newcommand{\ttilde}[1]{\widetilde{#1}}

\theoremstyle{plain}
\newtheorem{theointro}{Theorem}

\title[Braid groups of projective complex reflection groups]{Proof of Shvartsman's conjecture on braid groups of projective complex reflection groups}
\author{Owen Garnier}
\address{Departamento de \'Algebra e Instituto de Matem\'aticas de la Universidad de Sevilla, Spain.}
\email{owen.garnier@math.cnrs.fr}
\date{\today}

\begin{document}

\begin{abstract}
The purpose of this note is to prove a conjecture of Shvartsman relating a complex projective reflection group with the quotient of a suitable complex braid group by its center. Shvartsman originally proved this result  in the case of real projective reflection groups, and we extend it to all complex projective reflection groups. Our study also allows us to correct a result of Broué, Malle, Rouquier on projective reflection groups.
\end{abstract}

\maketitle
%\tableofcontents

\section*{Introduction}
Let $V$ be a finite dimensional complex vector space, and let $\P(V)$ be the associated projective space. The image of a subset $X\subset V$ in $\P(V)$ will be denoted by $\hhat{X}$. For $x\in V$, the line spanned by $x$ will be denoted by $[x]$. For $k$ a positive integer, the group of $k$-th complex roots of unity will be denoted by $\mu_k$. 

Consider a complex reflection group $W\subset \GL(V)$. In order to define the braid group $B$ attached to $W$, one first considers the subset $X$ of $V$ consisting of elements with a trivial stabilizer under the action of $W$, and then take the fundamental group of the quotient space $X/W$. 

Similarly, if $G\subset \PGL(V)$ is a \emph{projective reflection group} (i.e. a finite group which is generated by images of reflections in $\PGL(V)$), then we can consider the subset $\hhat{X}$ of $\P(V)$ consisting of elements with a trivial stabilizer under the action of $G$, and then take the fundamental group of the quotient space $\hhat{X}/G$. This gives a reasonable definition of a braid group attached to the projective reflection group $G$.

In his article \cite{shvartsman}, Shvartsman describes the braid group attached to a real projective reflection group as the quotient of a spherical Artin group by its center. More precisely, if $G\subset \PGL(V)$ is a real projective reflection group, then there is a maximal real reflection group $W\subset \GL(V)$ which is a lift of $G$. The braid group attached to $G$ is then isomorphic to the quotient $A/Z(A)$, where $A$ denotes the Artin group attached to $W$ \cite[Theorem A]{shvartsman}. The goal of Shvartsman is to use this realization in order to compute the possible orders of torsion elements in $A/Z(A)$, which is done in \cite[Theorem B]{shvartsman}.

In \cite[Section 4]{shvartsman}, Shvartsman conjectures that his results could be adapted in the case of complex projective reflection groups, where the Artin groups should be replaced by a complex braid group. The second result \cite[Theorem B]{shvartsman} about the possible orders of torsion elements was generalized to complex braid groups in \cite[Theorem 12.4]{beskpi1} and in \cite[Proposition 8.2]{regularbraids}. We give in this note an extension to all complex projective reflection groups of \cite[Theorem A]{shvartsman}:

\begin{theointro}[(Theorem \ref{theo:computation_projective})]\label{theo:shvartsman_intro}
Let $G\subset \PGL(V)$ be an irreducible projective reflection group, and let $\hhat{X}$ be the subset of $\P(V)$ consisting of elements with a trivial stabilizer under the action of $G$. Let also $W\subset \GL(V)$ be the maximal reflection group which is a lift of $G$. \\The fundamental group of $\hhat{X}/G$ is isomorphic to the quotient $B/Z(B)$, where $B$ denotes the complex braid group attached to $W$.
\end{theointro}

In Section \ref{sec:reminders_crg}, we give preliminary results about complex (projective) reflection groups. In particular, we prove that a projective complex reflection group admits a unique maximal lift in $\GL(V)$ which is a reflection group. In Section \ref{sec:main_results_mais_section_2_pas_la_sous-section}, we prove Theorem \ref{theo:shvartsman_intro} by introducing the \emph{enlarged braid group} attached to a complex reflection group $W$. This group is an extension of the more classical complex braid group attached to $W$, and we give some of its properties. Lastly, in Section \ref{sec:correction_bmr}, we give some corrections on a result of Broué, Malle, Rouquier in \cite{bmr}, which aims to relate braid groups of projective reflection groups with complex braid groups. As pointed out in \cite{dmm}, their result is false in general, and we give a complete description of the cases for which it holds (see Proposition \ref{prop:2.23bmr_reborn}).

\vspace{1em}\noindent\textbf{Acknowledgment.} We thank the Referee for their helpful comments, in particular for providing a concise and case free argument for Lemma \ref{lem:existence_releve_groupe_reflexion}, thus allowing us to considerably simplify the exposition of Section \ref{sec:reminders_crg}.

\section{Preliminaries}\label{sec:reminders_crg}
\subsection{Reminders and generalities on complex reflection groups}
In this section, we fix a finite dimensional complex vector space $V$, and we fix $n$ to be the dimension of $V$. We mostly follow \cite{lehrertaylor} for classical results on complex reflection groups.

Recall that a \emph{complex reflection group} $W$ is a finite subgroup of $\GL(V)$ which is generated by \emph{reflections}, that is finite order linear automorphisms of $V$ which pointwise fix some hyperplane.

A complex reflection group $W\subset \GL(V)$ is \emph{irreducible} if there are no $W$-invariant subspaces in $V$ apart from $\{0\}$ and $V$ itself. Every complex reflection group decomposes as a direct product of irreducible complex reflection groups and we will thus restrict our attention to irreducible groups from now on.  Irreducible complex reflection groups were classified by Shephard and Todd in \cite{shetod}, and we freely use the notation of \cite[Theorem 8.29]{lehrertaylor} regarding this classification. 

%In this classification, the groups $G(m,p,n)$ are exactly the \emph{imprimitive} groups. That is, if $W\subset \GL(V)$ is an irreducible complex reflection group such that $V$ decomposes as a direct sum of subspaces which are permuted by $W$, then $W$ has type $G(m,p,n)$ for some integers $m,p$. 

To an irreducible complex reflection group $W\subset \GL(V)$, one can attach the sequence $d_1\leqslant \ldots \leqslant d_n$ of its \emph{degrees}, and the sequence $d_n^*\geqslant \ldots \geqslant d_1^*=0$ of its \emph{codegrees} \cite[Proposition 3.25 and Definition 10.27]{lehrertaylor}. The degrees of $W$ are by definition the degrees of a \emph{system of basic invariants} of $W$ \cite[Theorem 3.20]{lehrertaylor}, that is, a family $(f_1,\ldots,f_n)$ of homogeneous elements of $S(V^*)$ which freely generate $S(V^*)^{W}$. Such a sequence always exists by the Chevalley-Shephard-Todd Theorem \cite[Theorem 3.20]{lehrertaylor}.  Moreover, a system of basic invariants $(f_1,\ldots,f_n)$ for $W$ induces a homeomorphism $V/W\simeq \C^n$, sending an orbit $W.v$ to the $n$-tuple $(f_1(v),\ldots,f_n(v))$ \cite[Proposition 9.3]{lehrertaylor}.

To a complex reflection group $W$ one can also attach the complement $X$ in $V$ of the union of the reflecting hyperplanes attached to the reflections of $W$. %A classical theorem of Steinberg \cite[Theorem 9.44]{lehrertaylor} ensures that the action of $W$ on $X$ is free and that we have a covering map $X\twoheadrightarrow X/W$. 
The \emph{braid group} $B(=B(W))$ of $W$ is then defined as the fundamental group of $X/W$, while the \emph{pure braid group} $P(=P(W))$ of $W$ is defined as the fundamental group of $X$ \cite[Section 2.B]{bmr}. The projection map $X\twoheadrightarrow X/W$ is a covering map by Steinberg's Theorem \cite[Theorem 9.44]{lehrertaylor}, and it induces a short exact sequence
\[1\to P\to B\to W\to 1.\]

Let $W\subset \GL(V)$ be a complex reflection group, and let $\zeta\in \C^*$. An element $g\in W$ is said to be $\zeta$-\emph{regular} if it admits a $\zeta$-eigenvector which lies in $X$. In other words, $g\in W$ is $\zeta$-regular if $V(g,\zeta)\cap X\neq \varnothing$, where $V(g,\zeta)$ denotes the $\zeta$-eigenspace of $g$. The eigenspace $V(g,\zeta)$ is then called a \emph{regular eigenspace} for $W$. An integer $k\geqslant 0$ is said to be \emph{regular} for $W$ if $W$ contains some $\zeta_k$-regular element, where $\zeta_k:=\exp(\frac{2i\pi}{k})$.

If $g\in W$ is $\zeta$-regular, then $g^k$ is $\zeta^k$ regular for all integer $k$. In particular, if $\zeta\in \C^*$ has order $m$, then $W$ contains a $\zeta$-regular element if and only if $m$ is regular for $W$. We also have that regular eigenspaces for $W$ can all be written as $V(g,\zeta_m)$, where $g$ is a $\zeta_m$-regular element for some integer $m$. 

An important criterion for regularity is given in \cite[Theorem 11.28]{lehrertaylor}, stating that a positive integer $k$ is regular for $W$ if and only if it divides as many degrees of $W$ as it does codegrees (both counted with multiplicity).

\subsection{Projective reflection groups}
In this section, we fix a finite dimensional complex vector space $V$ of dimension $n>1$.

We define a \emph{projective reflection group} as a finite subgroup of $\PGL(V)$ which is generated by images of reflections in $\PGL(V)$. Imitating the case of linear reflection group, we say that a projective reflection group $G$ is \emph{irreducible} if the only $G$-invariant subspace of $\P(V)$ is $\P(V)$ itself. %We say that $G$ is \emph{imprimitive} if there is a direct sum $V=V_1\oplus\cdots\oplus V_m$ with $m\geqslant 2$ such that the action of $G$ on $\P(V)$ permutes the subspaces $V_1,\ldots,V_m$ among themselves (we call $\{V_1,\ldots,V_m\}$ a \emph{system of imprimitivity} for $G$).

If $W\subset \GL(V)$ is an irreducible complex reflection group, then the image $\hhat{W}$ of $W$ in $\PGL(V)$ is a projective complex reflection group. The kernel of the projection map $W\twoheadrightarrow \hhat{W}$ is the set of scalar matrices lying in $W$. Since $W$ is irreducible, this set coincides with the center $Z(W)$ of $W$ by Schur's lemma, and we identify $\hhat{W}$ with the quotient group $W/Z(W)$. An elementary result is that every projective reflection group is the image in $\PGL(V)$ of some linear reflection group:

\begin{lem}\label{lem:existence_releve_groupe_reflexion}
Let $G\subset \PGL(V)$ be a projective reflection group. There is a complex reflection group $W\subset \GL(V)$ such that $\hhat{W}=G$. %Moreover, $W$ is imprimitive if and only if $G$ is primitive.
\end{lem}
\begin{proof}
Let $r\in \GL(V)$. For $\zeta\in \C^*$, the multiset of eigenvalues of $\zeta r$ is simply $\zeta$ times the multiset of eigenvalues of $r$. If some $\zeta r$ is a reflection in $\GL(V)$, then $r$ has at most two distinct eigenvalues $\lambda_1,\lambda_2$. Moreover, $n>1$ implies that $1$ is an eigenvalue of $\zeta r$, thus we have $\zeta\in \lambda_1^{-1},\lambda_2^{-1}$. In particular there is a finite number of elements of the form $\zeta r$ (i.e. lifts of the image of $r$ in $\PGL(V)$) which are reflections.

Let now $R$ denote the (finite) set of elements of $G$ which are images of reflections in $\GL(V)$. By the above argument, the set $\ttilde{R}$ of all reflections in $\GL(V)$ whose image in $\PGL(V)$ belongs to $R$ is finite. Let $W$ be the reflection group generated by $\ttilde{R}$. The image $\hhat{W}$ of $W$ in $\PGL(V)$ is equal to $G$ by construction, and it only remains to show that $W$ is finite.

Let $H$ be the subgroup of scalar matrices in $W$. The group $G$ is isomorphic to $W/H$. Since $G$ is finite, it is sufficient to show that $H$ is finite. Since $\ttilde{R}$ is a finite set of reflections, its image under the determinant lies in a finite group $\mu$ of roots of unity in $\C^*$. Since $W$ is generated by $\ttilde{R}$, the image of $W$ under the determinant also lies in $\mu$. Since the determinant of the scalar matrix $\zeta\Id$ is $\zeta^n$, we have $H\subset \{\zeta\Id~|~\zeta^n\in \mu\}$, which is finite.
\end{proof}

\begin{rem}
Let $W\subset \GL_1(\C)\simeq \C^*$ be a complex reflection group of rank 1. In this case $W$ is a cyclic group and the quotient $\hhat{W}$ is the trivial group acting on $\P^0(\C)=\{*\}$. This situation is not very rich and Theorem \ref{theo:shvartsman_intro} is immediate in this case. This is why we assume that $\dim V>1$ when considering projective reflection groups.
\end{rem}

In order to describe the action of $\hhat{W}$ on the projective space $\P(V)$, we first consider the inverse image $\ttilde{W}$ of $\hhat{W}$ in $\GL(V)$. We can describe the group $\ttilde{W}$ explicitly: Let $Z$ denote the center of $\GL(V)$, that is the subgroup of scalar multiples of the identity. The groups $W$ and $Z$ normalize each other, and we can consider the product group $ZW$ (not to be confused with the center $Z(W)$ of $W$). Since $Z$ is also the kernel of the natural morphism $\GL(V)\to \PGL(V)$, we have $ZW=\ttilde{W}$. This group contains both $Z$ and $W$ as normal subgroups, and we have :
\[\begin{cases}\ttilde{W}/Z\simeq W/(Z\cap W)=\hhat{W},\\ \ttilde{W}/W\simeq Z/(Z\cap W)\simeq \C^*/\mu_{|Z(W)|}\simeq \C^*.\end{cases}\]

During the proof of Lemma \ref{lem:existence_releve_groupe_reflexion}, we showed that the group $W_f$ generated by all the reflections lying in $\ttilde{W}$ is a complex reflection group. We call $W_f$ the \emph{full reflection group} associated to $W$. By construction, every complex reflection group $W'$ such that $\ttilde{W}=\ttilde{W'}$ is such that $W'\subset W_f$. Moreover, since $W$ is normal in $\ttilde{W}$, it is also normal in $W_f$.

The following result gives the type of $W_f$ depending on the type of $W$ when $W$ is irreducible.

\begin{prop}[\textbf{(Description of full reflection group)}]\label{prop:existence_full_group}
Let $W\subset \GL(V)$ be an irreducible complex reflection group and let $W_f$ be its full reflection group.
\begin{enumerate}
\item If $W$ is a primitive group and $\dim(V)=2$, then $W_f$ has type $G_7,G_{11},G_{19}$ depending on whether $W$ is a tetrahedral, octahedral or icosahedral group.
\item If $W$ is a primitive group and $\dim(V)>2$, then $W_f=W$, except when $W$ has type $G_{25}$, in which case $W_f$ has type $G_{26}$.
\item If $W$ has type $G(m,p,2)$, then $W_f$ has type $G(\frac{2m}{p\wedge 2},2,2)$, 
\item If $W$ has type $G(m,p,n)$ with $n>2$, then $W_f$ has type $G(m,p\wedge n,n)$.
\end{enumerate}
\end{prop}
\begin{proof}
Cases $(1)$ and $(2)$ are studied in \cite[Section 1.3 and 1.4]{shetod}, where the authors classify all the possible lifts of $\hhat{W}$.

Assume now that we are in cases $(3)$ or $(4)$, and let us denote by $W'$ the group we claim is equal to $W_f$. We have $W\subset W'$, and thus $Z(W)=Z(W')\cap W$ as both $W$ and $W'$ are irreducible. We can then see $\hhat{W}$ as a subgroup of $\hhat{W'}$. Since we have $[W:Z(W)]=[W':Z(W')]$, we deduce that $\hhat{W}=\hhat{W'}$. Now, the projection map $\GL(V)\to \PGL(V)$ induces a bijection between the subgroups of $\PGL(V)$ and the subgroups of $\GL(V)$ which contain $Z$. Since $\ttilde{W}$ and $\ttilde{W'}$ both contain $Z$, and since the image of $\ttilde{W}$ (resp. of $\ttilde{W'}$) in $\PGL(V)$ is equal to $\hhat{W}$ (resp. to $\hhat{W'}$), we have $\ttilde{W}=\ttilde{W'}$. Since $W'$ is also a complex reflection group, we have $W'\subset W_f$. 

$(3)$ The group $W$ has order $2m/(p\wedge 2)$, while $Z(W)$ has order $\frac{m}{p}(p\wedge 2)$. The quotient group $\hhat{W}$ then has order $\frac{2m}{p\wedge 2}$. Now, since $n=2$, any nontrivial element of $\hhat{W}$ admits exactly two lifts in $\GL(V)$ which are reflections. Thus $\ttilde{W}$ contains $\frac{4m}{p\wedge 2}-2$ reflections, as well as $W_f$. The number of reflections of $W'$ is the sum of its degrees minus 2, which is also equal to $\frac{4m}{p\wedge 2}-2$. Thus $W'$ contains all the reflections in $W_f$ and $W'=W_f$ as claimed.

$(4)$ The group $W=G(m,p,n)$ consists of monomial matrices, that is, matrices with exactly one nonzero entry by row and by column. The property of being monomial is invariant under multiplication by a scalar matrix, and thus $\ttilde{W}$ contains only monomial matrices, as well as $W_f$. We then have $W_f=G(m_f,p_f,n)$ for some integers $m_f,p_f$. Let $k\geqslant 1$. We can consider the reflection
\[M_k:=\matrix{0&\zeta_k&0\\\zeta_k^{-1}&0&0\\0&0&I_{n-2}}.\]
We have $M_m\in W\subset W_f$, which implies that $m$ divides $m_f$. Conversely, $M_{m_f}\in W_f$ implies that $\lambda M_{m_f}\in W$ for some $\lambda\in \C^*$. Since $n>2$, $\lambda$ is a coefficient of $\lambda M_{m_f}$, and thus we have $\lambda,\lambda \zeta_{m_f}\in \mu_{m}$. We deduce that $\zeta_{m_f}\in \mu_m$ and thus $m_f$ divides $m$. Since $W\subset W_f$ and since $\hhat{W}=\hhat{W_f}$, we have 
\[[W:Z(W)]=[W_f:Z(W_f)]\Leftrightarrow \frac{m}{p\wedge n}=\frac{m}{p_f\wedge n},\]
and $W_f$ ranges among the family of groups $\{G(m,p',n)~|~p'\wedge n=p\wedge n\}$. All the groups in this family are included in $W'$ and thus $W_f\subset W'$ as we wanted to show.
\end{proof}

\begin{rem}
Note that Proposition \ref{prop:existence_full_group} fails if $n=1$. Indeed, if $W$ is a complex reflection group of rank $1$, then $W$ is cyclic, and $\ttilde{W}=\C^*$ does not depend on $W$. This is another reason why we assume $n>1$ when considering projective reflection groups.
\end{rem}

We finish this section by proving that full reflection groups admit no regular hyperplanes. This result will be used in Section \ref{sec:main_results}.

\begin{lem}\label{lem:full_implies_no_regular_hyperplanes}Let $W\subset \GL(V)$ is a complex reflection group. If $W=W_f$ is its own full reflection group, then $W$ admits no regular hyperplanes in $V$.
\end{lem}
\begin{proof}
Let $\zeta\in \C^*$, and let $g\in W$ be an element whose $\zeta$-eigenspace is a hyperplane $H$. The element $\zeta^{-1}g\in \ttilde{W}$ pointwise fixes $H$: it is a reflection. We have $\zeta^{-1}g\in \ttilde{W}$ and thus $\zeta^{-1}g\in W_f=W$. In other words, $H$ is a reflecting hyperplane for $W_f=W$. In particular, $H$ cannot be a regular eigenspace for $W$.
\end{proof}

\section{Main results}\label{sec:main_results_mais_section_2_pas_la_sous-section}
In this section, we fix a finite dimensional complex vector space $V$ of dimension $n\geqslant 1$. We also fix an irreducible complex reflection group $W\subset \GL(V)$. We otherwise keep the notation from Section \ref{sec:reminders_crg}.

\subsection{Actions of $\ttilde{W}$, $\hhat{W}$ and $Z/Z(W)$}\label{sec:action_of_ttildew_hhatw_andzsurzw}
In order to prove Theorem \ref{theo:shvartsman_intro}, we need to study the action of $\hhat{W}$ on $\P(V)$. In this section, we relate this action with the action of $\ttilde{W}$ on $V$, and with the action of $\ttilde{W}/W$ on $V/W$. 

Since $Z$ is a normal subgroup of $\ttilde{W}$, there is an action of $\ttilde{W}/Z\simeq \hhat{W}$ on $V/Z$, which we can restrict to an action on $(V\setminus 0)/Z=\P(V)$. This action coincides with the natural action of $\hhat{W}$ on $\P(V)$. Similarly, since $W$ is a normal subgroup of $\ttilde{W}$, there is an action of $\ttilde{W}/W\simeq Z/Z(W)$ on $V/W$.

Note that the faithful action of $\hhat{W}$ on $\P(V)$ can also be seen as an action of $W$ on $\P(V)$ with kernel $Z(W)$. In particular, we will sometimes denote $\P(V)/W$ instead of $\P(V)/\hhat{W}$ to alleviate notation. Similarly, we can see the action of $Z/Z(W)$ on $V/W$ as an action of $Z\simeq \C^*$. We have the following commutative square of topological spaces:

% https://q.uiver.app/#q=WzAsNCxbMCwwLCJWXFxzZXRtaW51cyAwIl0sWzEsMCwiKFZcXHNldG1pbnVzIDApL1ciXSxbMSwxLCIoVlxcc2V0bWludXMgMCkvXFx0dGlsZGV7V30iXSxbMCwxLCJcXFAoVikiXSxbMCwxLCIvVyJdLFsxLDIsIi9aIl0sWzAsMywiL1oiLDJdLFszLDIsIi9XIiwyXSxbMCwyLCIvXFx0dGlsZGV7V30iXV0=
\begin{equation}\label{equ:commsquare}
\begin{tikzcd}
	{V\setminus 0} & {(V\setminus 0)/W} \\
	{\P(V)} & {(V\setminus 0)/\ttilde{W}}
	\arrow["{/W}", from=1-1, to=1-2]
	\arrow["{/Z}"', from=1-1, to=2-1]
	\arrow["{/\ttilde{W}}", from=1-1, to=2-2]
	\arrow["{/Z}", from=1-2, to=2-2]
	\arrow["{/W}"', from=2-1, to=2-2]
\end{tikzcd}\end{equation}
This square (or rather its restriction to a convenient subspace of $V$) will prove useful to us in the next section. Another relation between the actions of $\ttilde{W},\hhat{W}$ and $Z/Z(W)$ is the computation of nontrivial stabilizers:

\begin{lem}[\textbf{(Stabilizers)}]\label{lem:stabilizers}
Let $x\in V$ be different from $0$.
\begin{enumerate}[(a)]
\item The stabilizer of $x\in V$ under the action of $\ttilde{W}$ is given by
\[\Stab_{\ttilde{W}}(x)=\{\zeta w\in \ttilde{W}~|~x\in V(w,\zeta^{-1})\}.\]
It is nontrivial if and only if $x$ belongs to either a reflecting hyperplane of $W$ or to a proper regular eigenspace of $W$.
\item The stabilizer of $[x]$ under the action of $\hhat{W}$ is trivial if and only if $\Stab_{\ttilde{W}}(x)$ is trivial.
\item The stabilizer of $W.x$ under the action of $Z/Z(W)$ is trivial if and only if $\Stab_{\ttilde{W}}(x)\subset W$.
\end{enumerate}
\end{lem}

\begin{proof}
$(a)$ Let $\zeta w\in \ttilde{W}$. We have
\[(\zeta w).x=x\Leftrightarrow w.x=\zeta^{-1}x\Leftrightarrow x\in V(w,\zeta^{-1}).\]
%Note that, if $w=\omega\Id\in Z(W)$, we have $V(w,\omega)=V$, but we also have $\omega^{-1}w=\omega^{-1}\omega \Id=\Id$, thus central elements of $W$ do not give rise to nontrivial stabilizing elements in the action of $\ttilde{W}$.
Now, assume that $\Stab_{\ttilde{W}}(x)$ is nontrivial and that $x$ does not belong to a reflecting hyperplane of $W$. By assumption, there is some $w\in W$, and some $\zeta\in \C^*$ such that $x\in V(w,\zeta)$ and that $w\neq \zeta$ (otherwise $\zeta^{-1}w$ is trivial). Since $x$ does not belong to any reflecting hyperplane of $W$, $x\in V(w,\zeta)$ implies that $w$ is $\zeta$-regular. Then, $w\neq \zeta$ implies that $w$ is noncentral, and thus $V(w,\zeta)\subsetneq V$ is a proper regular eigenspace.

Conversely, assume that $x$ belongs to the reflecting hyperplane associated to a reflection $r\in W$. By definition, we have $x\in V(r,1)$, and thus $r\in \Stab_{\ttilde{W}}(x)$. Lastly, assume that $g$ is a $\zeta$-regular element such that $x\in V(g,\zeta)\subsetneq V$. We have $g\neq \zeta$ since $V(\zeta\Id,\zeta)=V$, and thus $1\neq \zeta^{-1}g\in \Stab_{\ttilde{W}}(x)$, which is nontrivial.

$(b)$ Let $wZ(W)\in \hhat{W}$. By definition, we have $wZ(W).[x]=[w.x]$ and thus
\begin{align*}
wZ(W)\in \Stab_{\hhat{W}}([x])&\Leftrightarrow [x]=[w.x]\\
&\Leftrightarrow \exists \zeta \in \C^*~|~w.x=\zeta^{-1} x.\\
&\Leftrightarrow \exists \zeta\in \C^*~|~(\zeta w).x=x\\
&\Leftrightarrow \exists w\zeta\in wZ\cap \Stab_{\ttilde{W}}(x).
\end{align*}
In particular, we obtain that $\Stab_{\hhat{W}}([x])$ is trivial if and only if $Z$ is the only $Z$-coset in $\ttilde{W}$ which intersects $\Stab_{\ttilde{W}}(x)$ nontrivially. This is equivalent to $\Stab_{\ttilde{W}}(x)\subset Z$. However, since $x$ is nonzero, $Z\cap \Stab_{\ttilde{W}}(x)$ is always equal to $\{1\}$, and $\Stab_{\ttilde{W}}(x)\subset Z$ is equivalent to $\Stab_{\ttilde{W}}(x)=\{1\}$.

$(c)$ Let $\lambda W\in \ttilde{W}/W$. By definition, we have $\lambda Z(W).W.x=W.(\lambda.x)$ and thus
\begin{align*}
\lambda W\in \Stab_{Z/Z(W)}(W.x)&\Leftrightarrow W.(\lambda x)=W.x\\
&\Leftrightarrow \exists w\in W~|~w\lambda x=x\\
&\Leftrightarrow \exists w\lambda \in \lambda W\cap \Stab_{\ttilde{W}}(x).
\end{align*}
In particular, we obtain that $\Stab_{Z/Z(W)}(W.x)$ is trivial if and only if $W$ is the only $W$-coset in $\ttilde{W}$ which intersects $\Stab_{\ttilde{W}}(x)$ nontrivially. This is equivalent to $\Stab_{\ttilde{W}}(x)\subset W$.
\end{proof}

\subsection{Enlarged complex braid groups}\label{sec:main_results} Now that we understand the stabilizers under the action of $\hhat{W}$ on $\P(V)$, we are ready to introduce the enlarged (pure) braid group of a complex reflection group.

\begin{definition}[\textbf{(Strongly regular vectors)}]
The set of \emph{strongly regular vectors} attached to $W$ is defined as 
\[X_S(=X_S(W)):=\{x\in V~|~\Stab_{\ttilde{W}}(x)=1\}.\]
\end{definition}

By Lemma \ref{lem:stabilizers}, $X_S$ is the complement in $V$ of the union of both the reflecting hyperplanes of $W$ and of its proper regular eigenspaces. We can now define the enlarged (pure) braid group as the fundamental group of $X_S/W$ (of $X_S$). We fix a basepoint $x_0$ in $X_S$ for the remainder of this section.

\begin{definition}[\textbf{(Enlarged braid group)}]
The \emph{enlarged braid group} attached to $W$ is defined as $B_S(=B_S(W)):=\pi_1(X_S/W,W.x_0)$ and the \emph{enlarged pure braid group} attached to $W$ is defined as $P_S(=P_S(W)):=\pi_1(X_S,x_0)$.
\end{definition}

The inclusion map $X_S/W\hookrightarrow X/W$ (resp $X_S\hookrightarrow X$) induces a morphism $B_S\to B$ (resp. $P\to P_S$). The following result gives some general information on these two morphisms:

\begin{prop}\label{prop:enlarged_braid_groups_is_an_extension_of_braid_group}\hfill
\begin{enumerate}[(a)]
\item The morphism $B_S\to B$ (resp. $P_S\to P$) induced by the inclusion $X_S/W\hookrightarrow X/W$ (resp. $X_S\hookrightarrow X$) is surjective. Moreover, if $W=W_f$, then it is an isomorphism.
\item In general, the group $B_S$ is the inverse image of $W$ under the natural projection $B(W_f)\twoheadrightarrow W_f$. It is a normal subgroup of index $[W_f:W]$ of $B(W_f)$. In particular, $B_S$ is torsion free. 
\end{enumerate}
\end{prop}
\begin{proof}
$(a)$ Among the proper regular eigenspaces of $W$, we distinguish between the regular hyperplanes and the regular eigenspaces of complex codimension $\geqslant 2$. Let $X_S'$ be the space obtained from $X$ by removing the regular hyperplanes of $W$. We have $X_S\subset X_S'\subset X$.

The space $X_S'/W$ (resp. $X_S'$) is obtained from $X/W$ (resp. from $X$) by removing an algebraic hypersurface. Since $X/W$ and $X$ are themselves complements of algebraic hypersurface in an affine space, \cite[Proposition A1]{bmr} gives that the inclusion $X_S'/W\hookrightarrow X/W$ (resp. $X_S'\hookrightarrow X$) induces a surjective morphism between the associated fundamental groups. Then, the space $X_S/W$ (resp. $X_S$) is obtained from the smooth complex manifold $X_S'/W$ (resp. $X_S'$) by removing a subvariety of complex codimension $\geqslant 2$. By \cite[Theorem X.2.3]{godbillon}, the inclusion $X_S/W\hookrightarrow X_S'/W$ (resp. $X_S\hookrightarrow X_S$) induces an isomorphism between the associated fundamental groups.

Now, if $W=W_f$, then $W$ admits no regular hyperplanes by Lemma \ref{lem:full_implies_no_regular_hyperplanes}. In this case, we have $X_S'=X$ and the second part of the above argument gives the isomorphism $B_S\simeq B$.

$(b)$ Consider $W_f$ the full reflection group attached to $W$. Since $X_S(W)$ depends only on $\ttilde{W}$, and since $\ttilde{W}=\ttilde{W_f}$, we have $X_S=X_S(W)=X_S(W_f)$. Since $W$ has finite index in $W_f$, the covering map $X_S\twoheadrightarrow X_S/W_f$ factors through the covering map $X_S\twoheadrightarrow X_S/W$ into a covering map $X_S/W\twoheadrightarrow X_S/W_f$ and $B_S$ is a finite-index subgroup of $B_S(W_f)$. Moreover, we know that $W$ is a normal subgroup of $W_f$, and thus the covering map $X_S/W\twoheadrightarrow X_S/W_f$ is the quotient by the action of $W_f/W$ on $X_S/W$. This covering map induces a short exact sequence
\[1\to B_S\to B_S(W_f)\to W/W_f \to 1\]
which gives the desired result.
\end{proof}

Now, restricting the commutative square (\ref{equ:commsquare}) to $X_S$ yields a commutative square:
\[\begin{tikzcd}
	X_S & {X_S/W} \\
	{\hhat{X_S}} & {\hhat{X_S}/W}
	\arrow[from=1-1, to=1-2]
	\arrow[from=1-1, to=2-1]
	\arrow[from=1-2, to=2-2]
	\arrow[from=2-1, to=2-2]
\end{tikzcd}\]
in which all the maps are fibrations:
\begin{itemize}
\item The fiber bundle $V\twoheadrightarrow \P(V)$ restricts to a fiber bundle $X_S\to \hhat{X_S}$. We obtain a short exact sequence
\[1\to \pi_1([x_0]\setminus 0,x_0)\to P_S\to \pi_1(\hhat{X_S},[x_0])\to 1.\]
In other words, the natural morphism $P_S\to \pi_1(\hhat{X_S},[x_0])$ is surjective and its kernel is generated by the (homotopy class in $X_S$) of the loop $t\mapsto \exp(2i\pi t)x_0$, which we denote by $\pi_S$.
\item The covering map $X\twoheadrightarrow X/W$ restricts to a covering map $X_S\twoheadrightarrow X_S/W$. We obtain a short exact sequence
\[1\to P_S\to B_S\to W\to 1.\] 
\item The action of $\hhat{W}$ on $\hhat{X_S}$ is free by Lemma \ref{lem:stabilizers}. Since $\hhat{W}$ is finite, it acts properly on $\hhat{X_S}$, which is locally compact. The projection $\hhat{X_S}\twoheadrightarrow \hhat{X_S}/W$ is then a covering map, and we have a short exact sequence
\[1\to \pi_1(\hhat{X_S},[x_0])\to \pi_1(\hhat{X_S}/W,W.[x_0])\to \hhat{W}\to 1.\]
\item The action of $\C^*\simeq Z/Z(W)$ on $X_S/W$ is free by Lemma \ref{lem:stabilizers}. It is then a free and proper Lie group action of $\C^*$ on $X_S/W$, which is a smooth manifold (as an open subset of $V/W\simeq\C^n$). By \cite[Theorem 21.10]{lee_manif}, the quotient map $p:X_S/W\to(X_S/W)/\C^*\simeq \hhat{X_S}/W$ is a principal $\C^*$-bundle and we have a short exact sequence
\[1\to \pi_1(p^{-1}(W.x_0),W.x_0)\to B_S\to \pi_1(\hhat{X_S}/W,W.[x_0])\to 1.\]
In other words, the natural morphism $B_S\to \pi_1(\hhat{X_S}/W,W.[x_0])$ is surjective and its kernel is generated by the (homotopy class in $X_S/W$) of the loop $t\mapsto \exp(2i\pi t).(W.x_0)$. By construction of the action of $\C^*\simeq Z/Z(W)$ on $V/W$ (see Section \ref{sec:action_of_ttildew_hhatw_andzsurzw}), this last loop is actually the loop $t\mapsto W.\left(e^{\frac{2i\pi t}{|Z(W)|}}x_0\right)$. We denote its homotopy class by $\beta_S$.
\end{itemize}

By construction, $\beta_S^{|Z(W)|}$ is the image in $B_S$ of $\pi_S$, and the image of $\beta_S$ in $W$ is $e^{\frac{2i\pi}{|Z(W)|}}\Id$, which generates $Z(W)$. These elements allow us to describe the center of enlarged braid groups: 

The center of complex braid groups was studied in \cite{bmr}, \cite{beskpi1} and \cite{dmm}. In the space $X$, the path $t\mapsto \exp(\frac{2i\pi t}{|Z(W)|})x$ induces a well defined element $\beta\in Z(B)$. The element $\beta^{|Z(W)|}=\pi\in P$ is represented by the path $t\mapsto \exp(2i\pi t)x$. The main results of \cite{dmm} state that $Z(B)$ (resp. $Z(P)$) is cyclic and generated by $\beta$ (resp. by $\pi$) when $W$ is irreducible. Moreover, if $U\subset B$ is a finite index subgroup, then $Z(U)\subset Z(B)$. Using their results, along with Theorem \ref{theo:main_diagram}, we are able to describe the center of enlarged braid groups.

\begin{cor}[\textbf{(Center of enlarged braid groups)}]
The center of $B_S$ is infinite cyclic and generated by $\beta_S$. If $U\subset B_S$ is a finite index subgroup, then $Z(U)\subset Z(B_S)$. In particular, the center of $P_S$ is infinite cyclic and generated by $\pi_S$.
\end{cor}
\begin{proof}
First, if $W=W_f$ has no regular hyperplane, then the natural morphism $B_S\to B$ (resp $P_S\to P$) is an isomorphism by Proposition \ref{prop:enlarged_braid_groups_is_an_extension_of_braid_group}. Since this morphism sends $\beta_S$ to $\beta$ (resp. $\pi_S$ to $\pi$), the result is precisely \cite[Theorem 1.2 and Theorem 1.4]{dmm}.

In the general case, we consider $W_f$ the full reflection group attached to $W$. We have $X_S=X_S(W)=X_S(W_f)$. By Proposition \ref{prop:enlarged_braid_groups_is_an_extension_of_braid_group}, $B_S$ is a finite index subgroup of $B(W_f)$. This shows directly that if $U\subset B_S$ has finite index, then $Z(U)\subset Z(B_S)$.

From the first part of the proof, we know that the center of $B_S(W_f)$ is infinite cyclic and generated by $\beta_S(W_f)$ and that the center of $B_S$ is infinite cyclic and generated by the smallest power of $\beta_S(W_f)$ which belongs to it. 

For a positive integer $k$, $\beta_S(W_f)^k$ is the homotopy class of the image in $X_S/W_f$ of the path in $X_S$ given by
\[\gamma_k:t\mapsto e^{\frac{2ik\pi t}{|Z(W_f)|}}x_0.\]
We have $\beta_S(W_f)^k\in B_S$ if and only if the endpoint $e^{\frac{2i k\pi}{|Z(W_f)|}}$ of $\gamma_k$ lies in $W.x_0$. Now, since $x_0\in X_S(W)$, having $\lambda x \in W.x$ for some $\lambda\in \C^*$ is possible only if $\lambda\Id\in Z(W)$, in other words if $\lambda$ is a power of $e^{\frac{2i\pi t}{|Z(W)|}}$. Thus, the endpoint of $\gamma_k$ lies in $W.x_0$ if and only if $k$ is a multiple of $[Z(W_f):Z(W)]$. The center of $B_S$ is then generated by the homotopy class of the image in $X_S/W$ of the path 
\[\gamma_{[Z(W_f):Z(W)]}:t\mapsto e^{\frac{2i\pi t}{|Z(W)|}}x_0,\]
which is precisely $\beta_S(W)$.

Lastly, $Z(P_S)$ is generated by the smallest power of $\beta_S(W)$ which lies in $P_S$. In other words, $Z(P_S)$ is generated by the homotopy class in $X_S$ of the path $\gamma_{k[Z(W_f):Z(W)]}$, where $k$ is the smallest integer such that the endpoint of $\gamma_{k[Z(W_f):Z(W)]}$ is $x_0$. Since the endpoint of $\gamma_{k[Z(W_f):Z(W)]}$ is $e^{\frac{2i\pi k}{|Z(W)|}}$, the smallest such integer is $k=|Z(W)|$, and $Z(P_S)$ is generated by the homotopy class in $X_S$ of the path $\gamma_{|Z(W_f)|}$, which is $\pi_S$ by definition.
\end{proof}

This corollary was the last argument needed, along with our commutative square of fibrations, to show the following result:

\begin{theo}\label{theo:main_diagram}Let $V$ be a finite dimensional complex vector space, and let $W\subset \GL(V)$ be an irreducible complex reflection group. Let also $x_0\in X_S$ be a basepoint. We have a commutative diagram, where all short sequences are exact:
% https://q.uiver.app/#q=WzAsOSxbMCwwLCJcXGxhbmdsZSBcXHBpX1NcXHJhbmdsZSJdLFsxLDAsIlxcbGFuZ2xlIFxcYmV0YV9TXFxyYW5nbGUiXSxbMiwwLCJaKFcpIl0sWzIsMSwiVyJdLFsxLDEsIkJfUyhXKSJdLFsxLDIsIlxccGlfMShcXGhoYXR7WH0vVyxXLlt4XzBdKSJdLFsyLDIsIlxcaGhhdHtXfSJdLFswLDIsIlxccGlfMShcXGhoYXR7WH0sW3hfMF0pIl0sWzAsMSwiUF9TKFcpIl0sWzAsMSwiIiwwLHsic3R5bGUiOnsidGFpbCI6eyJuYW1lIjoiaG9vayIsInNpZGUiOiJ0b3AifX19XSxbMSwyLCIiLDAseyJzdHlsZSI6eyJoZWFkIjp7Im5hbWUiOiJlcGkifX19XSxbMiwzLCIiLDAseyJzdHlsZSI6eyJ0YWlsIjp7Im5hbWUiOiJob29rIiwic2lkZSI6InRvcCJ9fX1dLFs0LDUsIiIsMCx7InN0eWxlIjp7ImhlYWQiOnsibmFtZSI6ImVwaSJ9fX1dLFszLDYsIiIsMCx7InN0eWxlIjp7ImhlYWQiOnsibmFtZSI6ImVwaSJ9fX1dLFs1LDYsIiIsMSx7InN0eWxlIjp7ImhlYWQiOnsibmFtZSI6ImVwaSJ9fX1dLFs3LDUsIiIsMSx7InN0eWxlIjp7InRhaWwiOnsibmFtZSI6Imhvb2siLCJzaWRlIjoidG9wIn19fV0sWzgsNywiIiwxLHsic3R5bGUiOnsiaGVhZCI6eyJuYW1lIjoiZXBpIn19fV0sWzAsOCwiIiwxLHsic3R5bGUiOnsidGFpbCI6eyJuYW1lIjoiaG9vayIsInNpZGUiOiJ0b3AifX19XSxbOCw0LCIiLDEseyJzdHlsZSI6eyJ0YWlsIjp7Im5hbWUiOiJob29rIiwic2lkZSI6InRvcCJ9fX1dLFs0LDMsIiIsMSx7InN0eWxlIjp7ImhlYWQiOnsibmFtZSI6ImVwaSJ9fX1dLFsxLDQsIiIsMCx7InN0eWxlIjp7InRhaWwiOnsibmFtZSI6Imhvb2siLCJzaWRlIjoidG9wIn19fV1d
\[\begin{tikzcd}
	{\langle \pi_S\rangle} & {\langle \beta_S\rangle} & {Z(W)} \\
	{P_S} & {B_S} & W \\
	{\pi_1(\hhat{X_S},[x_0])} & {\pi_1(\hhat{X_S}/W,W.[x_0])} & {\hhat{W}}
	\arrow[hook, from=1-1, to=1-2]
	\arrow[hook, from=1-1, to=2-1]
	\arrow[two heads, from=1-2, to=1-3]
	\arrow[hook, from=1-2, to=2-2]
	\arrow[hook, from=1-3, to=2-3]
	\arrow[hook, from=2-1, to=2-2]
	\arrow[two heads, from=2-1, to=3-1]
	\arrow[two heads, from=2-2, to=2-3]
	\arrow[two heads, from=2-2, to=3-2]
	\arrow[two heads, from=2-3, to=3-3]
	\arrow[hook, from=3-1, to=3-2]
	\arrow[two heads, from=3-2, to=3-3]
\end{tikzcd}\]
In particular, the fundamental group of $\hhat{X_S}/W$ (resp. of $\hhat{X_S}$) is isomorphic to $B_S/Z(B_S)$ (resp. to $P_S/Z(P_S)$).
\end{theo}

This theorem allows us to complete the proof of Theorem \ref{theo:shvartsman_intro};

\begin{theo}[\textbf{(Computation of projective braid groups)}]\label{theo:computation_projective}
Let $V$ be a finite dimensional complex vector space, and let $G\subset \PGL(V)$ be a nontrivial irreducible projective reflection group. Let also $W\subset \GL(V)$ be the maximal reflection group such that $\hhat{W}=G$, and let $\hhat{X}:=\{[x]\in \P(V)~|~\Stab_G([x])=1\}$. The fundamental group of $\hhat{X}/G$ is isomorphic to $B/Z(B)$.
\end{theo}
\begin{proof}
By construction, we have that $W=W_f$ is its own full group. By Proposition \ref{prop:enlarged_braid_groups_is_an_extension_of_braid_group}, the morphism $B_S(W)\to B(W)$ induced by the inclusion $X_S/W\hookrightarrow X/W$ is an isomorphism. Now, by Lemma \ref{lem:stabilizers}, we have $\hhat{X}=\hhat{X_S}$, thus $\hhat{X}/G=\hhat{X_S}/\hhat{W}$, and the fundamental group of this space is isomorphic to $B_S/\langle \beta_S\rangle \simeq B/\langle \beta\rangle=B/Z(B)$ by Theorem \ref{theo:main_diagram}.\end{proof}

Theorem \ref{theo:main_diagram} only exists at the level of enlarged braid groups, and cannot be extended to braid groups in every case (see Section \ref{sec:correction_bmr} for more details). However, it can be extended in the case where all regular elements of $W$ are central. 

\begin{cor}\label{cor:x=x_s_then_trivial}
If all regular elements of $W$ are central, then we have a commutative diagram, where all short sequences are exact:
\[\begin{tikzcd}
	{\langle \pi\rangle} & {\langle \beta\rangle} & {Z(W)} \\
	{P} & {B} & W \\
	{\pi_1(\hhat{X},[x_0])} & {\pi_1(\hhat{X}/W,W.[x_0])} & {\hhat{W}}
	\arrow[hook, from=1-1, to=1-2]
	\arrow[hook, from=1-1, to=2-1]
	\arrow[two heads, from=1-2, to=1-3]
	\arrow[hook, from=1-2, to=2-2]
	\arrow[hook, from=1-3, to=2-3]
	\arrow[hook, from=2-1, to=2-2]
	\arrow[two heads, from=2-1, to=3-1]
	\arrow[two heads, from=2-2, to=2-3]
	\arrow[two heads, from=2-2, to=3-2]
	\arrow[two heads, from=2-3, to=3-3]
	\arrow[hook, from=3-1, to=3-2]
	\arrow[two heads, from=3-2, to=3-3]
\end{tikzcd}\]
In particular, the fundamental group of $\hhat{X}/W$ (resp. of $\hhat{X}$) is isomorphic to $B/Z(B)$ (resp. to $P/Z(P)$).
\end{cor}
\begin{proof}
The only statement which is not directly implied by Theorem \ref{theo:main_diagram} is the fact that $\langle \beta\rangle=Z(B)$ (resp. $\langle \pi\rangle=Z(P)$) , which is known for all irreducible complex braid groups \cite[Theorem 1.1 and Theorem 1.2]{dmm}.
\end{proof}

\begin{rem}There are very few cases in which all regular elements of $W$ are central. Since the cardinality of the centralizer of a $\zeta_k$-regular element of $W$ is the product of the degrees of $W$ which are divisible by $k$, we have that a regular element is central if and only if its center divides all the degrees of $W$. Since the gcd of the degrees of $W$ is the cardinality of $Z(W)$, we obtain that all regular elements of $W$ are central if and only if all regular numbers for $W$ divide $Z(W)$.

If the highest degree of $W$ is regular, then all regular elements of $W$ are central if and only if all the degrees of $W$ are equal. The only irreducible groups satisfying this are $G_{7},G_{11},G_{19}$, along with $G(2d,2,2)$ for $d>2$. 

Assume now that the highest degree of $W$ is not regular. If $W$ belongs to the infinite series, then we must have in particular $W=G(m,p,n)$ with $\frac{m}{p}>1$ (otherwise the highest degree is regular). In this case, regular numbers are exactly the divisors of $\frac{mn}{p}$. All regular elements of $W$ are then central if and only if $\frac{mn}{p}$ divides $|Z(W)|=\frac{m}{p}(p\wedge n)$, which is equivalent to saying that $n$ divides $p$. Lastly, the only exceptional group $W$ for which the highest degree is not regular is $G_{15}$. By \cite[Theorem 11.28]{lehrertaylor}, all regular elements in $G_{15}$ are central.

We obtained that all regular elements of $W$ are central if and only if $W$ has type either $G_{7},G_{11},G_{15}$, $G_{19}$ or $G(m,p,n)$ with $p|n$.
\end{rem}

\subsection{Correction of a result of Broué, Malle, Rouquier}\label{sec:correction_bmr}
The natural commutative square
% https://q.uiver.app/#q=WzAsNCxbMCwwLCJYIl0sWzEsMCwiWC9XIl0sWzEsMSwiXFxoaGF0e1h9L1xcaGhhdHtXfSJdLFswLDEsIlxcaGhhdHtYfSJdLFswLDFdLFsxLDJdLFswLDNdLFszLDJdXQ==
\[\begin{tikzcd}
	X & {X/W} \\
	{\hhat{X}} & {\hhat{X}/\hhat{W}}
	\arrow[from=1-1, to=1-2]
	\arrow[from=1-1, to=2-1]
	\arrow[from=1-2, to=2-2]
	\arrow[from=2-1, to=2-2]
\end{tikzcd}\]
Induces a commutative diagram of groups
% https://q.uiver.app/#q=WzAsOSxbMCwxLCJQKFcpIl0sWzEsMSwiQihXKSJdLFsxLDIsIlxccGlfMShcXGhoYXR7WH0vXFxoaGF0e1d9LFxcaGhhdHtXfS5beF0pIl0sWzAsMiwiXFxwaV8xKFxcaGhhdHtYfSxbeF0pIl0sWzAsMCwiXFxsYW5nbGUgXFxwaVxccmFuZ2xlIl0sWzEsMCwiXFxsYW5nbGUgXFxiZXRhXFxyYW5nbGUiXSxbMiwxLCJXIl0sWzIsMCwiWihXKSJdLFsyLDIsIlxcaGhhdHtXfSJdLFswLDFdLFsxLDJdLFswLDNdLFszLDJdLFs0LDBdLFs0LDVdLFs1LDFdLFsxLDZdLFs1LDddLFs3LDZdLFs2LDhdXQ==
\begin{equation}\label{eq:basic_diagram}
\begin{tikzcd}
	{\langle \pi\rangle} & {\langle \beta\rangle} & {Z(W)} \\
	{P} & {B} & W \\
	{\pi_1(\hhat{X},[x])} & {\pi_1(\hhat{X}/\hhat{W},\hhat{W}.[x])} & {\hhat{W}}
	\arrow[from=1-1, to=1-2]
	\arrow[from=1-1, to=2-1]
	\arrow[from=1-2, to=1-3]
	\arrow[from=1-2, to=2-2]
	\arrow[from=1-3, to=2-3]
	\arrow[from=2-1, to=2-2]
	\arrow[from=2-1, to=3-1]
	\arrow[from=2-2, to=2-3]
	\arrow[from=2-2, to=3-2]
	\arrow[from=2-3, to=3-3]
	\arrow[from=3-1, to=3-2]
\end{tikzcd}\end{equation}

The second row of this diagram is a short exact sequence (see Section \ref{sec:reminders_crg}). The first row is a short exact sequence by \cite[Theorem 1.3]{dmm}. The first column is a short exact sequence, induced by the $\C^*$-bundle $X\twoheadrightarrow \hhat{X}$ (restriction of the $\C^*$-bundle $V\setminus 0\twoheadrightarrow \P(V)$). The third column is a short exact sequence by construction.

Contrary to the projection map $X_S/W\twoheadrightarrow \hhat{X_S}/\hhat{W}$ studied in Section \ref{sec:main_results}, the projection map $X/W\twoheadrightarrow \hhat{X}/\hhat{W}$ may not be a fibration. In particular, we cannot use the same argument as in the proof of Theorem \ref{theo:main_diagram} to obtain that the second column of Diagram \eqref{eq:basic_diagram} is a short exact sequence. However, we still have the following partial result:

\begin{prop}\label{prop:surjectivity_natural_projection}
The morphism $B\to \pi_1(\hhat{X}/\hhat{W},\hhat{W}.[x])$ induced by the projection map $X/W\to \hhat{X}/\hhat{W}$ is surjective. Furthermore, if $b\in B$ admits a nontrivial power in $Z(B)$, then the image of $b$ in $\pi_1(\hhat{X}/\hhat{W},\hhat{W}.[x])$ is trivial.
\end{prop}
\begin{proof}
First, $\hhat{X}$ is algebraic open subset of $\P(V)$. In particular it is an irreducible normal variety. Since the quotient of a normal variety by a finite group is again normal, the quotient $\hhat{X}/\hhat{W}$ is again an irreducible normal variety.  Now, consider the open set $X_S\subset X$ studied in Section \ref{sec:main_results}. 
For $x\in X_S$, we have a commutative diagram  
% https://q.uiver.app/#q=WzAsNCxbMCwwLCJCX1MoVykiXSxbMSwwLCJYL1ciXSxbMCwxLCJcXHBpXzEoXFxoaGF0e1hfU30vXFxoaGF0e1d9LFxcaGhhdHtXfS5beF0pIl0sWzEsMSwiXFxwaV8xKFxcaGhhdHtYfS9cXGhoYXR7V30sXFxoaGF0e1d9Llt4XSkiXSxbMCwxXSxbMCwyXSxbMSwzXSxbMiwzXV0=
\begin{equation}\label{eq:diag_blabla}\begin{tikzcd}
	{B_S} & {B} \\
	{\pi_1(\hhat{X_S}/\hhat{W},\hhat{W}.[x])} & {\pi_1(\hhat{X}/\hhat{W},\hhat{W}.[x])}
	\arrow[from=1-1, to=1-2]
	\arrow[from=1-1, to=2-1]
	\arrow[from=1-2, to=2-2]
	\arrow[from=2-1, to=2-2]
\end{tikzcd}\end{equation}

We are going to show that all the morphisms in this diagram are surjective. Let us recall that, if $Y$ is an irreducible complex normal algebraic variety, and if $U\subset Y$ is a nonempty (algebraic) open subset of $Y$, then the natural morphism from the fundamental group of $U$ to that of $Y$ is surjective (see \cite[Theorem 2.1]{fundamental_group_normal_varieties} and the references there).

\begin{itemize}
\item By Proposition \ref{prop:enlarged_braid_groups_is_an_extension_of_braid_group}, the morphism $B_S\to B$ is surjective.
\item Since $\hhat{X}/\hhat{W}$ is a normal variety, and since $\hhat{X_S}/\hhat{W}$ is an algebraic open subset of $\hhat{X}/\hhat{W}$, the morphism $\pi_1(\hhat{X_S}/\hhat{W},\hhat{W}.[x])\to \pi_1(\hhat{X}/\hhat{W},\hhat{W}.[x])$ is surjective. 
\item The morphism $B_S\to \pi_1(\hhat{X_S}/\hhat{W},\hhat{W}.[x])$ was proven to be surjective in Theorem \ref{theo:main_diagram}.
\end{itemize}
Since three out of four arrows in Diagram \eqref{eq:diag_blabla} are surjective, then so is the fourth one, which is what we wanted to show.

Lastly, let $d$ be a regular number for $W$, let $g\in W$ be a $d$-regular element, and let $x\in V(g,\zeta_d)\cap X$. Fixing $x$ as a basepoint, we can consider the path $\gamma:[0,1]\to X$ sending $t$ to $\exp(\frac{2i\pi t}{d} )x$. The endpoint of $\gamma$ is $\zeta_d x=g.x\in W.x$. Thus $\gamma$ induces a well defined element $\ttilde{g}$ of $B$ which is a lift of $g$ under the projection map $B\twoheadrightarrow W$. By construction, $\ttilde{g}$ belongs to the kernel of the morphism $B\to \pi_1(\hhat{X}/\hhat{W},\hhat{W}.[x])$ since the image of the path $\gamma$ remains in $\ttilde{W}.x$ at all times. 

Now, the element $\ttilde{g}$ is a $d$-th root of the full-twist $\pi$. Moreover, by \cite[Theorem 1.2 and Proposition 8.1]{regularbraids}, an element $b\in B$ admits a central power if and only if it is conjugate to a power of an element of the form $\ttilde{g}$ (for some regular element $g\in W$). Thus, any element of $B$ admitting a central power belongs to the kernel of the morphism $B\twoheadrightarrow \pi_1(\hhat{X}/\hhat{W},\hhat{W}.[x])$.
%Lastly, $\beta$ is the homotopy class of the loop $t\mapsto W.(e^{2i\pi t/|Z(W)|}x)$. Since the image of this loop in $\hhat{X}/\hhat{W}$ is constant, the image of $\beta$ in $\pi_1(\hhat{X}/\hhat{W},\hhat{W}.[x])$ is trivial.
\end{proof}

Note that the above result applies to $\beta$, which admits itself as a central power in $B$. In particular, the composition of the second column of Diagram \eqref{eq:basic_diagram} is trivial, even if the sequence is not exact. 

Now, it is claimed in  \cite[Proposition 2.23]{bmr} that Diagram \eqref{eq:basic_diagram} can always be completed by a morphism $\pi_1(\hhat{X}/\hhat{W},\hhat{W}.[x])\to \hhat{W}$ into a commutative diagram in which every row and every column is a short exact sequence. It is pointed out in \cite{dmm} that this result is false in general. We give a complete description of the cases in which it holds:

\begin{prop}\label{prop:2.23bmr_reborn}
The following statements are equivalent:
\begin{enumerate}[(i)]
\item All regular elements in $W$ are central.
\item The second column in Diagram \eqref{eq:basic_diagram} is a short exact sequence. 
\item There is a morphism $\pi_1(\hhat{X}/\hhat{W},\hhat{W}.[x])\to \hhat{W}$ which completes Diagram \eqref{eq:basic_diagram} into a commutative diagram.
\end{enumerate}
Moreover, if these statements hold, then all the rows and columns in the completed diagram are short exact sequences. 
\end{prop}
\begin{proof}
Let $K\subset B$ denote the kernel of the natural morphism $B\to \pi_1(\hhat{X}/\hhat{W},\hhat{W}.[x])$. By Proposition \ref{prop:surjectivity_natural_projection}, the group $\pi_1(\hhat{X}/\hhat{W},\hhat{W}.[x])$ is isomorphic to the quotient $B/K$. Point $(iii)$ is then equivalent to stating that the kernel of the morphism $B\to \hhat{W}$ contains $K$. The kernel of the morphism $B\to \hhat{W}$ is the preimage of $Z(W)$ under the morphism $B\twoheadrightarrow W$: it is the subgroup of $B$ generated by $P$ and $\beta$.

$(i)\Rightarrow (ii)$ is proven in Corollary \ref{cor:x=x_s_then_trivial}.

$(ii)\Rightarrow (iii)$: If the second column of Diagram \eqref{eq:basic_diagram} is a short exact sequence, then $K=\langle \beta\rangle\subset \langle P,\beta\rangle$, which proves $(iii)$.

$(iii)\Rightarrow (i)$: we prove the contrapositive. Assume that $W$ admits a noncentral regular element $g$. By \cite[Theorem 1.2]{regularbraids}, there is some root of $\pi$ in $B$ which is a lift of $g$ under the projection map $B\twoheadrightarrow W$. By Proposition \ref{prop:surjectivity_natural_projection}, $\rho$ belongs to the kernel of the projection map $B\twoheadrightarrow \pi_1(\hhat{X}/\hhat{W},\hhat{W}.[x])$. However, since $g$ is not central in $W$, the image of $\rho$ in $\hhat{W}$ is not trivial. The kernel $K$ of the morphism $B\to \pi_1(\hhat{X}/\hhat{W},\hhat{W}.[x_0])$ is then not included in $\langle P,\beta\rangle$ and $(iii)$ is false.
Lastly, we saw in Corollary \ref{cor:x=x_s_then_trivial} that, if all regular elements in $W$ are central, then we can complete Diagram \eqref{eq:basic_diagram} into a commutative diagram in which all rows and all columns are short exact sequences.
\end{proof}

\begin{rem}In the general case, Proposition \ref{prop:surjectivity_natural_projection} only gives that we have a surjective morphism $B/K\twoheadrightarrow \pi_1(\hhat{X}/\hhat{W},\hhat{W}.[x])$, where $K$ is the subgroup of $B$ generated by the elements in $B$ which admit a central power. We do not know whether this morphism is an isomorphism in the general case. However, preliminary computations seem to imply that in most cases, $B/K$ is either trivial, or a cyclic group, which would in turn imply that $\pi_1(\hhat{X}/\hhat{W},\hhat{W}.[x])$ is either trivial or cyclic, at least in a number of cases. We chose not to include those rather long computations here, as they only give partial results.  
\end{rem}

\printbibliography
\end{document}